\documentclass[10pt]{article}

\usepackage[T2A]{fontenc}
\usepackage{amsfonts}
\usepackage{amsmath}
\usepackage{amssymb}
\usepackage{comment}
\usepackage[dvips]{graphicx}
\usepackage{subfigure}

\newtheorem{theorem}{Theorem}[section]
\newtheorem{lemma}[theorem]{Lemma}
\newtheorem{corollary}[theorem]{Corollary}

\newtheorem{prop}[theorem]{Proposition}
\newtheorem{remark}[theorem]{Remark}

\newcommand{\set}[1]{\left \{ #1 \right \}}                     
\newcommand{\setst}[2]{\left \{ #1 \mid #2 \right \}}           
\newcommand{\abs}[1]{\left| #1 \right|} 

\providecommand{\R}{\mathbb{R}}

\renewcommand{\tilde}{\widetilde}
\renewcommand{\hat}{\widehat}
\renewcommand{\bar}{\overline}

\newcommand{\calB}{\mathcal{B}}
\newcommand{\calF}{\mathcal{F}}

\newcommand{\deltain}{\delta^{\rm in}}
\newcommand{\deltaout}{\delta^{\rm out}}
	
\newenvironment{proof}{\par\noindent%
{\bf Proof.\par\nopagebreak}}{\unskip\nobreak\enskip$\square$\par\bigskip}
\newenvironment{proofof}[1]{\medskip\par\noindent%
{\bf Proof of #1.\par\nopagebreak}}{\unskip\nobreak\enskip$\square$\par\bigskip}

\newcommand{\reffig}[1]{Fig.~\ref{fig:#1}}           
\newcommand{\refth}[1]{Theorem~\ref{th:#1}}          
\newcommand{\reflm}[1]{Lemma~\ref{lm:#1}}            
\newcommand{\refsec}[1]{Section~\ref{sec:#1}}        
\newcommand{\refrem}[1]{Remark~\ref{rem:#1}}         

\makeindex

\begin{document}

\title
{
	Acyclic Bidirected and Skew-Symmetric Graphs:\\
	Algorithms and Structure
}

\author
{
	Maxim A. Babenko
	\thanks
	{
		Dept. of Mechanics and Mathematics,
		Moscow State University, Vorob'yovy Gory, 119899 Moscow,
		Russia, \textsl{email}: mab@shade.msu.ru.
		Supported by RFBR grants 03-01-00475 and NSh 358.2003.1.
	}
}

\maketitle

\begin{abstract}
	\emph{Bidirected graphs} (a sort of nonstandard graphs introduced
by Edmonds and Johnson) provide a natural generalization to the
notions of directed and undirected graphs. By a \emph{weakly
(node- or edge-) acyclic} bidirected graph we mean such graph
having no (node- or edge-) simple cycles.
We call a bidirected graph \emph{strongly acyclic} if it has
no cycles (even non-simple). Unlike the case of standard graphs,
a bidirected graph may be weakly acyclic but still have non-simple cycles.

Testing a given bidirected graph for weak acyclicity
is a challenging combinatorial problem, which also has a number
of applications (e.g. checking a perfect matching in a general
graph for uniqueness). We present (generalizing results of Gabow,
Kaplan, and Tarjan) a modification of the depth-first search
algorithm that checks (in linear time) if a given bidirected graph
is weakly acyclic (in case of negative answer a simple cycle is
constructed).

Our results are best described in terms of \emph{skew-symmetric graphs}
(the latter give another, somewhat more convenient graph language
which is essentially equivalent to the language of bidirected graphs).

We also give structural results for the class of weakly acyclic bidirected
and skew-symmetric graphs explaining how one can construct any such graph
starting from strongly acyclic instances and, vice versa, how one can
decompose a weakly acyclic graph into strongly acyclic ``parts''. Finally,
we extend acyclicity test to build (in linear time) such a
decomposition.

\end{abstract}	



\section{Introduction}
\label{sec:intro}

The notion of \emph{bidirected graphs} was introduced by
Edmonds and Johnson~\cite{EJ-70} in connection with one important class of integer linear programs
generalizing problems on flows and matchings; for
a survey, see also~\cite{law-76,sch-03}.

Recall that in a \emph{bidirected} graph $G$ three types of edges are allowed:
(i)~a~usual directed edge, or an \emph{arc}, that leaves one node and enters another one;
(ii)~an~edge \emph{from both} of its ends; or
(iii)~an~edge \emph{to both} of its ends.

When both ends of edge coincide, the edge becomes a loop.

In what follows we use notation $V_G$ (resp. $E_G$) to denote the set of
nodes (resp. edges) of an undirected or bidirected graph~$G$.
When $G$ is directed we speak of arcs rather than edges and write $A_G$ in place
of $E_G$.

A \emph{walk} in a bidirected graph $G$ is an alternating sequence
$P = (s = v_0, e_1, v_1, \ldots, e_k, v_k = t)$ of nodes and edges such that
each edge $e_i$ connects nodes $v_{i-1}$ and $v_i$, and for
$i = 1, \ldots, k-1$, the edges $e_i,e_{i+1}$ form a \emph{transit pair}
at $v_i$, which means that one of $e_i,e_{i+1}$ enters and the
other leaves~$v_i$. Note that $e_1$ may enter $s$ and $e_k$ may
leave $t$; nevertheless, we refer to $P$ as a walk from $s$ to $t$,
or an $s$--$t$ \emph{walk}. $P$ is a \emph{cycle} if $v_0=v_k$ and
the pair $e_1,e_k$ is transit at $v_0$; a cycle is usually
considered up to cyclic shifts. Observe that an $s$--$s$ walk
is not necessarily a cycle.

If $v_i \ne v_j$ for all $1 \le i < j < k$ and $1 < i < j \le k$,
then walk $P$ is called \emph{node-simple} (note that the endpoints of a
node-simple walk need not be distinct). A walk is called \emph{edge-simple}
if all its edges are different.

Let $X$ be an arbitrary subset of nodes of $G$.
One can modify $G$ as follows: for each node $v\in X$ and each edge
$e$ incident with $v$, reverse the direction of $e$ at $v$.
This transformation preserves the set of walks in $G$ and
thus does not change the graph in essence.
We call two bidirected graphs $G_1, G_2$ \emph{equivalent}
if one can obtain $G_2$ from $G_1$ by applying a number
of described transformations.
 
\begin{figure}[tb]
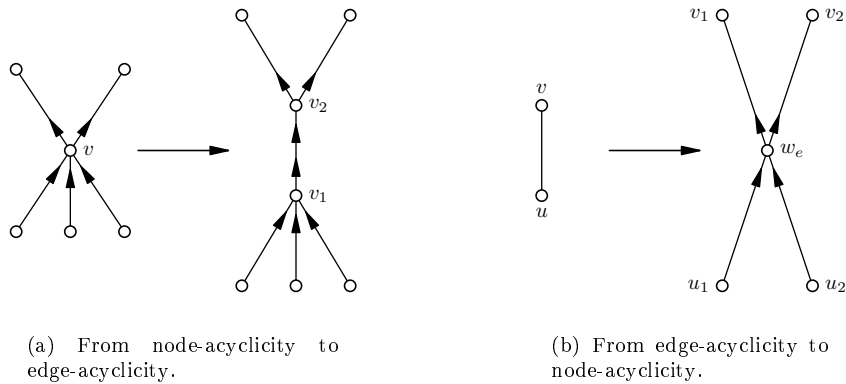

	\centering
	\subfigure[From node-acyclicity to edge-acyclicity.]{
		\includegraphics{pics/reductions.1}%
		\label{fig:node_to_edge}
	}
	\hspace{2cm}%
	\subfigure[From edge-acyclicity to node-acyclicity.]{
		\includegraphics{pics/reductions.2}%
		\label{fig:edge_to_node}
	}
	\caption{Reductions between the notions of weak node- and edge-acyclicity.}
\end{figure}

A bidirected graph is called \emph{weakly (node- or edge-) acyclic} if it has no
(\mbox{node-} or \mbox{edge-}) simple cycles. These two notions of acyclicity
are closely connected.
Given a bidirected graph~$G$ one can do the following:
(i) replace each node~$v \in V_G$ by a pair of nodes $v_1$, $v_2$;
(ii) for each node~$v \in V_G$ add an edge leaving $v_1$ and entering $v_2$;
(iii) for each edge $e \in E_G$ connecting nodes $u, v \in V_G$ add
an edge connecting $u_i$ and $v_j$, where $i = 1$ if $e$ enters $u$; $i = 2$ otherwise;
similarly for $j$ and $v$.
This procedure yields a weakly edge-acyclic graph
iff the original graph is weakly node-acyclic (see~\reffig{node_to_edge}).
The converse reduction from edge-acyclicity to node-acyclicity
is also possible:
(i)~replace each node~$v \in V_G$ by a pair of nodes $v_1$, $v_2$;
(ii)~for each edge $e \in E_G$ connecting nodes $u, v \in V_G$ add a node~$w_e$ and four edges
connecting $u_i$, $v_i$ with $w_e$ ($i = 1, 2$);
edges $u_iw_e$ should enter~$w_e$;
edges $w_ev_i$ should leave~$w_e$;
the directions of these edges at $u_i$ (resp.~$v_i$) should coincide
with the direction of $e$ at $u$ (resp.~$v$)
(see~\reffig{edge_to_node}).

In what follows we shall only study the notion
of weak edge-acyclicity. Hence, we drop the prefix ``edge'' for brevity
when speaking of weakly acyclic graphs.
If a bidirected graph has no (even non-simple) cycles we call it
\emph{strongly acyclic}. 

One possible application of weak acyclicity testing is described in \cite{GKT-99}.
Let $G$ be an undirected graph and $M$ be a \emph{perfect matching} in $G$
(that is, a set of edges such that: (i) no two edges in $M$ share a common node;
(ii) for each node~$v$ there is a matching edge incident with $v$).
The problem is to check if $M$ is the unique perfect matching in $G$.
To this aim we transform $G$ into the bidirected graph~$\bar G$ by assigning
directions to edges as follows: every edge $e \in M$ leaves both
its endpoints, every edge $e \in E_G \setminus M$
enters both its endpoints. One easily checks that the definition of matching implies
that every edge-simple cycle in $\bar G$ is also node-simple. Moreover,
each such simple cycle in $\bar G$ gives rise to an \emph{alternating circuit} in $G$ with respect to $M$
(a circuit of even length consisting of an alternating sequence of edges belonging to $M$
and $E_G \setminus M$). And conversely, every alternating circuit in $G$ with respect to $M$
generates a node-simple cycle in $\bar G$. It is well known (see \cite{LP-86})
that $M$ is unique iff there is no alternating circuit with respect to it.
Hence, the required reduction follows.

\section{Skew-Symmetric Graphs}
\label{sec:ss}

This section contains terminology and some basic facts concerning
skew-sym\-met\-ric graphs and explains the correspondence between these and
bidirected graphs. For a more detailed survey on skew-symmetric graphs,
see, e.g., ~\cite{tut-67,GK-96,GK-04,BK-05}.

A \emph{skew-symmetric graph} is a digraph~$G$ endowed with
two bijections $\sigma_V, \sigma_A$ such that: $\sigma_V$ is
an involution on the nodes (i.e., $\sigma_V(v)\ne v$ and
$\sigma_V(\sigma_V(v)) = v$ for each node~$v$), $\sigma_A$ is an
involution on the arcs, and for each arc $a$ from $u$ to $v$,
$\sigma_A(a)$ is an arc from $\sigma_V(v)$ to $\sigma_V(u)$. For
brevity, we combine the mappings $\sigma_V, \sigma_A$ into one mapping
$\sigma$ on $V_G \cup A_G$ and call $\sigma$ the \emph{symmetry} (rather
than skew-symmetry) of $G$.
For a node (arc) $x$, its symmetric node (arc) $\sigma(x)$ is also
called the \emph{mate} of $x$, and we will often use notation with
primes for mates, denoting $\sigma(x)$ by $x'$.

Observe that if $G$ contains an arc $e$ from a node $v$
to its mate $v'$, then $e'$ is also an arc from $v$ to $v'$ (so the
number of arcs of $G$ from $v$ to $v'$ is even and these parallel
arcs are partitioned into pairs of mates).

By a path (circuit) in $G$ we mean a node-simple
directed walk (cycle), unless explicitly stated otherwise.
The symmetry $\sigma$ is extended in a natural way to walks, cycles, paths,
circuits, and other objects in $G$. In particular, two walks or
cycles are symmetric to each other if the elements of one of them
are symmetric to those of the other and go in the reverse order:
for a walk $P = (v_0, a_1, v_1, \ldots, a_k, v_k)$, the symmetric
walk~$\sigma(P)$ is $(v'_k, a'_k, v'_{k-1}, \ldots, a'_1, v'_0)$.
One easily shows that $G$ cannot contain self-symmetric circuits
(cf.~\cite{GK-04}). We call a set of nodes~$X$ \emph{self-symmetric}
if~$X' = X$.

Following terminology in~\cite{GK-96}, an arc-simple walk in $G$ is
called \emph{regular} if it contains no pair of symmetric arcs
(while symmetric nodes in it are allowed). Hence, we may speak of
regular paths and regular circuits.

Next we explain the correspondence between skew-symmetric and
bidirected graphs (cf.~\cite[Sec.~2]{GK-04}, \cite{BK-05}). For sets $X, A, B$, we
use notation $X = A \sqcup B$ when $X = A \cup B$ and
$A \cap B= \emptyset$. Given a skew-symmetric
graph~$G$, choose an arbitrary partition $\pi=\set{V_1, V_2}$ of
$V_G$ such that $V_2$ is symmetric to $V_1$. Then $G$ and $\pi$ determine
the bidirected graph $\bar G$ with node set $V_1$ whose edges
correspond to the pairs of symmetric arcs in $G$. More precisely,
arc mates $a,a'$ of $G$ generate one edge $e$ of $\bar G$ connecting nodes
$u, v \in V_1$ such that: (i)~$e$~goes from $u$ to $v$ if one of $a, a'$
goes from $u$ to $v$ (and the other goes from $v'$ to $u'$ in
$V_2$); (ii)~$e$~leaves both $u,v$ if one of $a,a'$ goes from $u$
to $v'$ (and the other from $v$ to $u'$); (iii)~$e$~enters both
$u, v$ if one of $a,a'$ goes from $u'$ to $v$ (and the other from
$v'$ to $u$). In particular, $e$ is a loop if $a, a'$ connect a pair
of symmetric nodes.

\begin{figure}[tb]
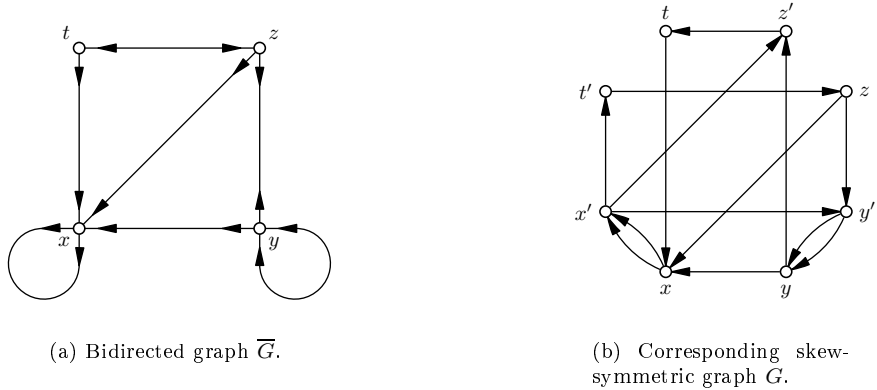

	\centering
	\subfigure[Bidirected graph~$\bar G$.]{
		\includegraphics{pics/examples.1}%
	}
	\hspace{3cm}%
	\subfigure[Corresponding skew-sym\-met\-ric graph~$G$.]{
		\includegraphics{pics/examples.2}%
	}
	\caption{Related bidirected and skew-symmetric graphs.}
	\label{fig:sk-bi}
\end{figure}

Conversely, a bidirected graph $\bar G$ with node set $\bar V$
determines a skew-sym\-met\-ric graph~$G$ with symmetry $\sigma$ as
follows. Take a copy $\sigma(v)$ of each element $v$ of $\bar V$,
forming the set $\bar V' := \setst{\sigma(v)}{v\in \bar V}$. Now set
$V_G := \bar V \sqcup \bar V'$. For each edge $e$ of $\bar G$ connecting nodes
$u$ and $v$, assign two ``symmetric'' arcs $a,a'$ in $G$ so as to satisfy
(i)-(iii) above (where $u' = \sigma(u)$ and $v' = \sigma(v)$). An
example is depicted in Fig.~\ref{fig:sk-bi}.

\begin{remark}
	A bidirected graph generates one skew-symmetric graph,
	while a skew-symmetric graph generates a
	number of bidirected ones, depending on the partition $\pi$ of $V$
	that we choose in the first construction. The latter bidirected graphs
	are equivalent.
\end{remark}

Also there is a correspondence between walks in $\bar G$ and walks in $G$.
More precisely, let $\tau$ be the natural mapping
of $V_G \cup A_G$ to $V_{\bar G} \cup E_{\bar G}$ (obtained by identifying the
pairs of symmetric nodes and arcs). Each walk
$$
	P = (v_0, a_1, v_1, \ldots, a_k, v_k)
$$
in $G$ induces the sequence
$$
	\tau(P) := (\tau(v_0), \tau(a_1), \tau(v_1), \ldots, \tau(a_k), \tau(v_k))
$$
of nodes and edges in $\bar G$. One can easily check that $\tau(P)$ is a walk in $\bar G$
and $\tau(P) = \tau(P')$. Moreover, for any walk $\bar P$ in $\bar G$
there is exactly one preimage $\tau^{-1}(\bar P)$ in $G$.

\medskip
Let us call a skew-symmetric graph \emph{strongly acyclic} if it has no directed cycles.
Each cycle in $\bar G$ generates a cycle in $G$ and vice versa.
To obtain a similar result for the notion of weak acyclicity in bidirected graphs,
suppose $\bar G$ is not weakly acyclic and consider an edge-simple cycle~$\bar C$ in $\bar G$ having the
smallest number of edges. Then $\bar C$ generates a cycles $C$ in $G$
(as described above). Cycles $C$, $C'$ are circuits since otherwise one can shortcut
them and obtain (by applying $\tau$) a shorter edge-simple cycle in $\bar G$.
Moreover, $C$ and $C'$ are regular (or, equivalently, arc-disjoint).
Indeed, suppose $C$ contains both arcs~$a$ and $a'$ for some $a \in A_G$.
Hence $\bar C$ traverses the edge $\tau(a)$ at least twice, contradicting the assumption.
Conversely, let $C$ be a regular circuit in $G$. Trivially $\bar C := \tau(C)$ is an edge-simple
cycle in $\bar G$. These observations motivate the following definition:
we call a skew-symmetric graph \emph{weakly acyclic} if is has no regular circuits.

The following proposition summarizes our observations.
\begin{prop}
\label{prop:acyclicity_red}
	$\bar G$ is strongly (resp. weakly) acyclic iff $G$ is strongly (resp. weakly) acyclic.
\end{prop}

For a given set of nodes~$X$ in a directed graph~$G$ we use notation $G[X]$ to denote the directed
subgraph induced by $X$. In case $G$ is skew-symmetric and $X' = X $ the symmetry on $G$
induces the symmetry on $G[X]$.

An easy part of our task is to describe the set of strongly acyclic skew-symmetric graphs.
The following theorem gives the complete characterization of such graphs.

\begin{theorem}
\label{th:strong_acycl_decomp}
	A skew-symmetric graph $G$ is strongly acyclic iff there exists a partition $Z \sqcup Z'$ of $V_G$,
	such that the induced (standard directed) subgraphs $G[Z]$, $G[Z']$ are acyclic
	and no arc goes from $Z$ to $Z'$.
\end{theorem}

In terms of bidirected graphs \refth{strong_acycl_decomp} means the following:
\begin{corollary}
	A bidirected graph $G$ is strongly acyclic iff $G$ is equivalent to a
	bidirected graph that only has directed edges forming an acyclic graph
	and edges leaving both endpoints.
\end{corollary}

\section{Separators and Decompositions}
\label{sec:separators}

In this section we try to answer the following question:
given a skew-symmetric weakly acyclic graph what kind of a natural
certificate can be given to prove the absence of regular circuits
(or, equivalently, regular cycles) in it? 

\begin{figure}[tb]
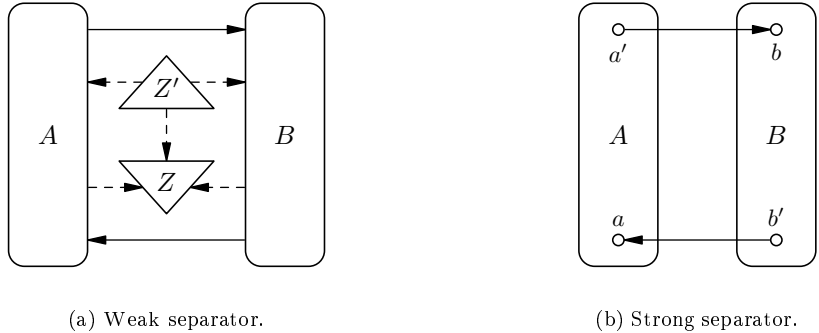

	\centering
	\subfigure[Weak separator.]{
		\includegraphics{pics/separator.1}%
		\label{fig:weak_separator}
	}
	\hspace{3cm}%
	\subfigure[Strong separator.]{
		\includegraphics{pics/separator.2}%
		\label{fig:strong_separator}
	}
	\caption{
		Separators. Solid arcs should occur exactly once,
		dashed arcs may occur arbitrary number of times (including zero)
	}
\end{figure}

Our first answer is as follows. Let $G$ be a skew-symmetric graph.
Suppose $V_G$ is partitioned into four sets $A, B, Z, Z'$ such that:
(i)~$A$ and $B$ are self-symmetric and nonempty;
(ii)~exactly one pair of symmetric arcs connects $A$ and $B$;
(iii)~$G[A]$ and $G[B]$ are weakly acyclic;
(iv)~no arc leaves $Z$, no arc enters $Z'$.
If these properties are satisfied we call $(A, B, Z)$ a \emph{weak separator} for~$G$
(see~\reffig{weak_separator}).

\begin{theorem}
\label{th:weak_separator}
	Every weakly acyclic skew-symmetric graph~$G$ is either strongly acyclic
	or admits a weak separator $(A, B, Z)$. Conversely,
	if $(A, B, Z)$ is a weak separator for $G$, then $G$ is
	weakly acyclic.
\end{theorem}

Thus, given a weakly acyclic graph~$G$ one can apply \refth{weak_separator} to 
split~$V_G$ into four parts. The subgraphs $G[A]$, $G[B]$ are again weakly acyclic, so
we can apply the same argument to them, etc. This recursive process (which produces
two subgraphs on each steps) stops when current subgraph becomes strongly acyclic.
In such case, \refth{strong_acycl_decomp} provides us with the required certificate.

Motivated by this observation we introduce the notion of a \emph{weak acyclic decomposition}
of $G$. By this decomposition we mean a binary tree $D$ constructed as follows.
The nodes of $D$ correspond to self-symmetric subsets of $V_G$
(in what follows, we make no distinction between nodes in $D$ and these subsets).
The root of $D$ is the whole node set $V_G$. Any leaf~$X$ in $D$ is a self-symmetric subset
that induces a strongly acyclic subgraph~$G[X]$; we attach a partition $X = Z \sqcup Z'$ as
in \refth{strong_acycl_decomp} to~$X$.
Consider any non-leaf node $X$ in $D$. It induces
the subgraph~$G[X]$ that is not strongly acyclic. Applying \refth{weak_separator}
we get a partition of $X$ into subsets
$A, B, Z, Z'$ and attach it to $X$; the children of $X$ are defined to be $A$ and $B$.

Provided that a weak separator can be found in linear time,
the above-described procedure for building weak acyclic decomposition totally requires
$O(mn)$ time in worst case ($n := \abs{V_G}$, $m := \abs{A_G}$).
However, one can use depth-first search
to construct a weak decomposition in linear time, see~\refsec{algorithms}.
This improved algorithm has a number of interesting applications.
For example, it can serve as a part of a procedure
that finds a shortest regular path in a weakly acyclic skew-symmetric graph
under arbitrary arc lengths and runs in $O(m \log^2 n)$ time.
However, this problem is quite complicated and will be addressed in another paper.

\medskip

An appealing special case arises when we restrict our attention to the class
of strongly connected (in a standard sense) skew-symmetric graphs,
that is, graphs where each two nodes are connected by a 
(not necessarily regular) path.

We need to introduce two additional definitions.
Given a skew-symmetric graph $H$ and a node $s$ in it we call $H$ \emph{$s$-connected}
if every node in $H$ lies on a (not necessarily regular) $s$--$s'$ path.
Suppose the node set of a skew-symmetric graph~$G$ admits a partition $(A,B)$ such that:
(i) $A$ and $B$ are self-symmetric;
(ii) exactly one pair of symmetric arcs $\set{a'b, b'a}$ connects $A$ and $B$ ($a, a' \in A$, $b, b' \in B$);
(iii) $G[A]$ is weakly acyclic and $a$-connected, $G[B]$ is weakly acyclic and $b$-connected.
Then we call $(A,B)$ a \emph{strong separator} for~$G$ (see~\reffig{strong_separator} for an example).

A simple corollary of \refth{weak_separator} is the following:
                                     
\begin{theorem}
\label{th:strong_separator}
	A skew-symmetric graph~$B$ is strongly connected and weakly acyclic
	iff it admits a strong separator $(A, B)$.
\end{theorem}

Now we extend \refth{strong_acycl_decomp} to describe a decomposition of
an arbitrary weakly acyclic skew-symmetric graph in terms of
strongly connected components (hence, providing another answer to the question
posed at the beginning of the section).

\begin{figure}[tb]
	\centering
	\includegraphics{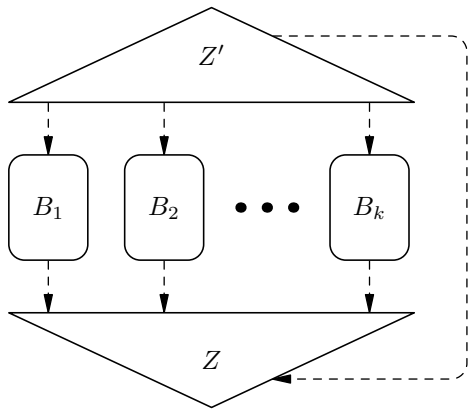}
	\caption{
		Decomposition of a weakly acyclic skew-symmetric graph~$G$.
		Dashed arcs may occur arbitrary number of times (including zero).
		Subgraphs $G[Z]$, $G[Z']$ are acyclic, subgraphs $G[B_i]$
		are strongly connected and weakly acyclic.
	}
	\label{fig:decomp}
\end{figure}
	
\begin{theorem}
\label{th:acycl_decomp}
	A skew-symmetric graph $G$  is weakly acyclic iff there exists a
	partition of $V_G$ into sets $Z, Z', B_1, \ldots, B_k$ such that:
	(i) (standard directed) subgraphs $G[Z]$, $G[Z']$ are acyclic;
	(ii) sets $B_i$ are self-symmetric, subgraphs $G[B_i]$ are strongly connected and weakly acyclic;
	(iii) no arc connects distinct sets $B_i$ and $B_j$;
	(iv) no arc leaves~$Z$, no arc enters~$Z'$.
\end{theorem}

An example of such decomposition is presented in \reffig{decomp}.
For $k = 0$ the decomposition in \refth{acycl_decomp} coincides
with such in \refth{strong_acycl_decomp}.

Consider an arbitrary weakly acyclic skew-symmetric graph~$G$.
Add auxiliary nodes $\set{s, s'}$ and arcs $\set{sv, v's'}$,
$v \in V_G \setminus \set{s,s'}$ thus making $G$ $s$-connected.
Similarly to its weak counterpart, a \emph{strong acyclic decomposition} of $G$ is a tree~$D$
constructed as follows. The nodes of $D$ correspond to self-symmetric subsets of $V_G$.
Each such subset~$A$ induces the $a$-connected graph~$G[A]$ for some $a \in A$.
The root of $D$ is the whole node set $V_G$. Consider a node $A$ of $D$. Applying \refth{acycl_decomp}
one gets a partition of $A$ into subsets $Z, Z', B_1, \ldots, B_k$
and attaches it to $A$. Each of $B_i$ is strongly connected and thus \refth{strong_separator} applies. Hence, we can
further decompose each of $B_i$ into $X_i \sqcup Y_i$ ($X_i' = X_i$, $Y_i' = Y_i$)
with the only pair of symmetric arcs $\set{x_i'y_i, y_i'x_i}$ ($x_i \in X_i$, $y_i \in Y_i$)
connecting $X_i$ and $Y_i$. The induced subgraphs $G[X_i]$ (resp. $G[Y_i]$)
are $x_i$-connected (resp. $y_i$-connected).
We define the children of $A$ to be $X_1, Y_1, \ldots, X_k, Y_k$.
Clearly, leaf nodes of $D$ correspond to certain strongly acyclic subgraphs.

The complexity of the described tree construction procedure is $O(mn)$
where $n := \abs{V_G}$, $m := \abs{A_G}$
(again, if linear-time algorithm for constructing separators is applied).
But unlike the case of weak decomposition, we are unaware of any faster algorithm. So building
strong decomposition in $o(mn)$ time is an open problem.

\section{Algorithms}
\label{sec:algorithms}

We need some additional notation. For a set of nodes~$X$ denote the set of arcs entering
(resp. leaving) $X$ by $\deltain(X)$ (resp. $\deltaout(X)$).
Denote the set of arcs having both endpoints in $X$ by $\gamma(X)$.

\begin{figure}[tb]
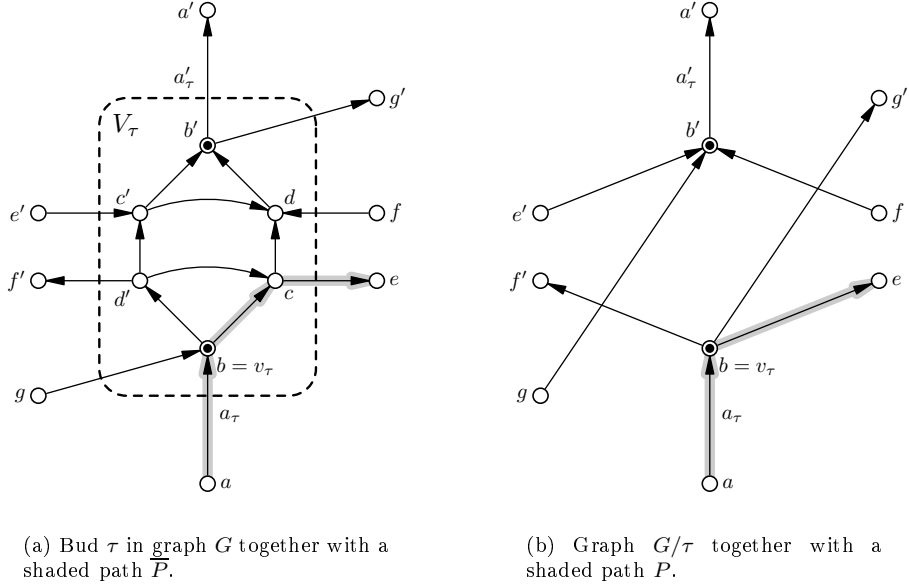

	\centering
	\subfigure[Bud~$\tau$ in graph~$G$ together with a shaded path~$\bar P$.]{
		\includegraphics{pics/trimming.1}%
		\label{fig:bud}
	}
	\hspace{1cm}%
	\subfigure[Graph~$G / \tau$ together with a shaded path $P$.]{
		\includegraphics{pics/trimming.2}%
	}
	\caption{Buds, trimming, and path restoration.
	Base and antibase nodes $b, b'$ are marked. Path $\bar P$ is a preimage of~$P$.}
	\label{fig:trimming}
\end{figure}

Let $V_\tau$ be a symmetric set of nodes in a skew-symmetric graph~$G$;
$a_\tau \in \deltain(V_\tau)$.
Let $v_\tau$ denote the head of $a_\tau$.
Suppose every node in $V_\tau$ is reachable from $v_\tau$ by a regular path in $G[V_\tau]$.
Then we call $\tau = (V_\tau, a_\tau)$ a \emph{bud}. (Note that our definition of bud
is weaker than the corresponding one in \cite{GK-96}.) The arc~$a_\tau$ (resp. node~$v_\tau$)
is called the \emph{base arc} (resp. \emph{base node}) of $\tau$,
arc~$a_\tau'$ (resp. node~$v_\tau'$) is called the \emph{antibase arc} (resp.
\emph{the antibase node}) of $\tau$. For an arbitrary bud~$\tau$ we denote its
set of nodes by $V_\tau$, base arc by $a_\tau$, and base node by $v_\tau$.
An example of bud is given in \reffig{bud}.

Consider an arbitrary bud~$\tau$ in a skew-symmetric graph~$G$.
By \emph{trimming~$\tau$} we mean the following transformation of $G$:
(i) all nodes in $V_\tau \setminus \set{v_\tau, v_\tau'}$ and arcs in $\gamma(V_\tau)$
are removed from $G$;
(ii) all arcs in $\deltain(V_\tau) \setminus \set{a_\tau}$ are transformed into
arcs entering $v_\tau'$ (the tails of such arcs are not changed);
(iii) all arcs in $\deltaout(V_\tau) \setminus \set{a_\tau'}$ are transformed into
arcs leaving $v_\tau$ (the heads of such arcs are not changed).
The resulting graph (which is obviously skew-symmetric) is denoted by $G / \tau$.
Thus, each arc of the original graph $G$ not belonging to $\gamma(V_\tau)$
has its \emph{image} in the trimmed graph $G / \tau$.
\reffig{trimming} gives an example of bud trimming.

Let $P$ be a regular path in $G / \tau$. One can lift this path to $G$ as follows:
if $P$ does not contain neither $a_\tau$, nor $a_\tau'$ leave~$P$ as it is.
Otherwise, consider the case when $P$ contains~$a_\tau$ (the symmetric case is analogous).
Split $P$ into two parts: the part $P_1$ from the beginning of $P$ to $v_\tau$ and
the part~$P_2$ from $v_\tau$ to the end of $P$. Let $a$ be the first arc of $P_2$. The arc~$a$
leaves $v_\tau$ in $G / \tau$ and thus corresponds to some arc $\bar a$ leaving $V_\tau$ in $G$
($\bar a \ne a_\tau'$). Let $u \in V_\tau$ be the tail of $a$ in $G$ and $Q$ be
a regular $v_\tau$--$u$ path in $G[V_\tau]$ (existence of $Q$ follows from definition of bud).
Consider the path $\bar P := P_1 \circ Q \circ P_2$ (here $U \circ V$ denotes the
path obtained by concatenating $U$ and $V$). One can easily show that $\bar P$ is regular.
We call $\bar P$ a \emph{preimage of~$P$} (under trimming $G$ by $\tau$).
Clearly, $\bar P$ is not unique.
An example of such path restoration is shown in \reffig{trimming}: the shaded path~$\bar P$
on the left picture corresponds to the shaded path~$P$ on the right picture.

Given a skew-symmetric graph~$G$ we check if it is weakly
acyclic as follows (we refer to this algorithm as \textsc{Acyclicity-Test}).
For technical reasons we require $G$ to obey the following
two properties:
\begin{enumerate}
	\item[(i)] \emph{Degree property}: for each node~$v$ in $G$ at most one arc enters~$v$
	or at most one arc leaves~$v$.  
	\item[(ii)] \emph{Loop property}: $G$ must not contain parallel arcs connecting
	symmetric nodes (these arcs correspond to loops in bidirected graphs).
\end{enumerate}	

Degree property implies that a regular walk in $G$ cannot contain a pair of
symmetric nodes (loosely speaking, the notions of node- and arc-regularity
coincide for $G$).

\begin{remark}
\label{rem:degree_reduction}
	Observe, that for a graph obtained by applying node- to edge-acyclicity
	reduction (as described~in \refsec{intro}) the degree property holds.
	Hence, to check an arbitrary graph for node-acyclicity 
	one may apply that reduction and invoke \textsc{Acyclicity-Test}.
	Similarly, to check a graph for edge-acyclicity we first
	reduce the problem to checking for node-acyclicity and then proceed as
	described earlier.
\end{remark}

Our algorithm adopts ideas from \cite{GKT-99} to the case of skew-symmetric graphs.
The algorithm is a variation of both depth-first-search (DFS) procedure (see~\cite{CLR-90}) and
regular reachability algorithm (see~\cite{GK-96}).
It has, however, two essential differences. Firstly,
unlike standard DFS, which is carried out in a static graph,
our algorithm changes~$G$ by trimming some buds.
Secondly, unlike regular reachability algorithm, we do not trim a bud as soon
as we discover it. Rather, trimming is postponed up to a moment when it can be
done ``safely''.

Degree and loop properties are preserved by trimmings. Indeed, consider
a bud~$\tau$ in a current graph~$H$. The node~$v_\tau$ has at least two
outgoing arcs (since there are two arc-disjoint $v_\tau$--$v_\tau'$ paths
in $H[V_\tau]$). Hence, $v_\tau$ has exactly one incoming arc (namely, $a_\tau$).
When $\tau$ is trimmed the in- and out-degrees can only change for $v_\tau$
and $v_\tau'$. For the node~$v_\tau$ (resp. $v_\tau'$) its
in- (resp. out-) degree remains~1, and thus degree property still holds.
Loop property is also maintained since trimming cannot produce parallel arcs
between base and antibase nodes.

Let $H$ be a current graph. Each pair of symmetric nodes in $G$ is mapped to a certain pair
of symmetric nodes in $H$. This mapping is defined by induction on the number
of trimmings performed so far. Initially this mapping is identity. When a bud $\tau$ is trimmed
and nodes $V_\tau \setminus \set{v_\tau,v_\tau'}$ are removed, the mapping
is changed so as to send the pairs of removed nodes to $\set{v_\tau,v_\tau'}$.
Given this mapping, we may also speak of the \emph{preimage}~$\bar X$ of any self-symmetric node
set~$X$ in $H$.

The algorithm recursively grows a directed forest~$F$.
At every moment this forest has no symmetric nodes
(or, equivalently, does not intersect the symmetric forest~$F'$).
Thus, every path in such forest is regular.
The algorithm assigns \emph{colors} to nodes.
There are five possible colors: white, gray, black, antigray, and antiblack.
White color assigned to $v$ means that $v$ is not yet discovered. Since the algorithm
processes nodes in pairs, if $v$ is white then so is $v'$. 
Other four colors also occur in pairs: if $v$ is gray then $v'$ is antigray,
if $v$ is black then $v'$ is antiblack (and vice versa).
All nodes outside both $F$ and $F'$ are white, nodes in $F$ are black or gray,
nodes in $F'$ are antiblack or antigray.

At any given moment the algorithm implicitly maintains a regular path starting from a root of~$F$.
As in usual DFS, this path can be extracted by examining the recursion stack.
The nodes on this path are gray, the symmetric nodes are antigray.
No other node is gray or antigray. Black color denotes nodes which are already completely
processed by the algorithm; the mates of such nodes are antiblack.

The core of the algorithm is the following recursive procedure.
It has two arguments~--- a node~$u$ and optionally an arc~$q$ entering~$u$
($q$ may be omitted when $u$ is a root node for a new tree in $F$).
Firstly, the procedure marks $u$ as gray and adds $u$ to $F$ (together with~$q$
if $q$ is given). Secondly, it scans all arcs leaving $u$. Let $a$ be such arc, $v$ be its head.
Several cases are possible (if no case applies, then $a$ is skipped and next
arc is fetched and examined):

\begin{enumerate}
	\item[(i)] \emph{Circuit case:}
	If $v$ is gray, then there exists a regular circuit in the current
	graph (it can be obtained by adding the arc~$a$ to the gray
	$v$--$u$ path in $F$). The procedure halts reporting
	the existence of a regular circuit in $G$ (which is constructed from~$C$
	in a postprocessing stage, see below).

	\item[(ii)] \emph{Recursion case:}
	If $v$ is white, the recursive call with parameters $(v, uv)$ is made.

	\item[(iii)] \emph{Trimming case:}
	If $v$ is antiblack, the procedure constructs a certain bud in the current graph
	and trims it as follows. We shall prove in the sequel that each time trimming case occurs
	the node~$v'$ is an ancestor of $u$ in $F$. Let $P$ denote the corresponding $u$--$v'$ path.
	Let $a_\tau$ be the (unique) arc of $F$ entering $u$ ($u$ has at least two outgoing arcs
	and hence cannot the a root of $F$, see below).
	Let $H$ denote the current graph.
	Finally, let $V_\tau$ be the union of node sets of $P$ and $P'$. One can easily show
	that $\tau = (V_\tau, a_\tau)$ is a bud in $H$ (buds
	formed by a pair of symmetric regular paths are called \emph{elementary} in \cite{GK-96}).
	The procedure trims $\tau$ and replaces $H$ by
	$H / \tau$. The forest~$F$ is updated by removing nodes in $V_\tau \setminus \set{u,u'}$
	and arcs in $\gamma(V_\tau)$. All other arcs of $F$ are replaced by their images under
	trimming by~$\tau$. Since $a_\tau$ belongs to $F$, it follows that
	the structure of forest is preserved. Note that trimming can produce new (previously unexisting)
	arcs leaving~$u$.
\end{enumerate}

When all arcs leaving $u$ are fetched and processed the procedure
marks $u$ as black, $u'$ as antiblack and exits.

\textsc{Acyclicity-Test} initially makes all nodes white.
Then, it looks for symmetric pairs of white nodes in $G$. Consider such pair $\set{v,v'}$
and assume, without loss of generality, that out-degree of $v$ is at most~1.
Invoke the above-described procedure at $v$ (passing no arc)
and proceed to the next pair.

If all recursive calls complete normally, we claim that the initial
graph is weakly acyclic. Otherwise, some recursive call halts
yielding a regular circuit~$C$ in a current graph.
During the postprocessing stage we consider the sequence of
the trimmed buds in the reverse order and undo the corresponding trimmings. Each time we undo
trimming of some bud $\tau$ we also replace $C$ by its preimage
(as described in \refsec{separators}). At each such step the regularity of $C$ is
preserved, thus at the end of postprocessing we obtain a regular circuit in the original
graph, as required.

The correctness of the algorithm will be proved in \refsec{correctness},
a linear-time implementation is given in \refsec{impl}.

\medskip

Now we address the problem of building a weak acyclic decomposition. 
We solve it by the algorithm \textsc{Decompose} which is
a modified version of \textsc{Acyclicity-Test}.

\begin{figure}[tb]
	\centering
	\includegraphics{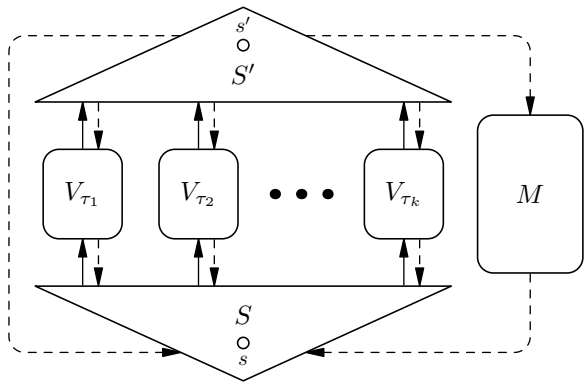}%
	\caption{
		A barrier. Solid arcs should occur exactly once,
		dashed arcs may occur arbitrary number of times (including zero).
	}
	\label{fig:barrier}
\end{figure}

Let $G$ be a skew-symmetric graph with a designated node $s$.
Suppose we are given a collection of buds $\tau_1, \ldots, \tau_k$ in $G$ together with
node sets $S$ and $M$. Additionally, suppose the following properties hold:
(i) collection $\set{S, S', M, V_{\tau_1}, \ldots, V_{\tau_k}}$ is a partition of~$V_G$
with $s \in S$;
(ii) no arc goes from $S$ to $S' \cup M$;
(iii) no arc connects distinct sets $V_{\tau_i}$ and $V_{\tau_j}$;
(iv) no arc connects $V_{\tau_i}$ and $M$;
(v) the arc~$e_{\tau_i}$ is the only one going from $S$ to $V_{\tau_i}$.
Then we call the tuple $\calB = (S, M; \tau_1, \ldots, \tau_k)$ 
an $s$--$s'$ \emph{barrier} (\cite{GK-96}, see~\reffig{barrier}
for an example).

Let us introduce one more weak acyclicity certificate
(which is needed for technical reasons) and show how to construct
a weak decomposition from it.
Let $\calB = (S, M; \tau_1, \ldots, \tau_k)$ be a barrier in~$G$.
Put $\tilde G := G / \tau_1 / \ldots / \tau_k$, $W := S \cup \set{v_{\tau_1}, \ldots, v_{\tau_k}}$.
We call $\calB$ \emph{acyclic} if the following conditions are satisfied:
(i) subgraphs $G[M], G[V_{\tau_1}], \ldots, G[V_{\tau_k}]$ are weakly acyclic.
(ii) the (standard directed) subgraph $\tilde G[W]$ is acyclic.

\begin{figure}[tb]
	\centering
	\includegraphics{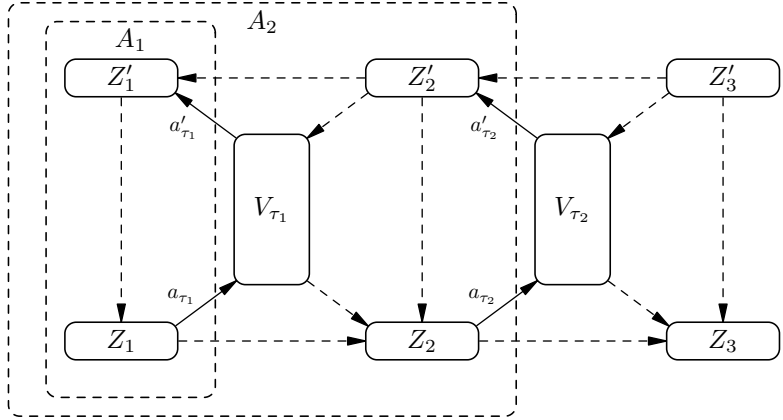}%
	\caption{
		Constructing a weak decomposition from an acyclic barrier.
		Solid arcs should occur exactly once, dashed arcs may occur
		arbitrary number of times (including zero).
		Not all possible dashed arcs are shown.
	}
	\label{fig:barrier_to_decomp}
\end{figure}

Suppose we are given an acyclic barrier~$\calB$ of $G$ with $M = \emptyset$.
Additionally, suppose weak acyclic decompositions of $G[V_{\tau_i}]$ are also given.
A weak acyclic decomposition of $G$ can be obtained as follows.
Consider the graph $\tilde G$ and the set~$W$ as in definition of an acyclic barrier.
Order the nodes in $W$ topologically: $W = \set{w_1, \ldots, w_n}$;
for $i > j$ no arc in $\tilde G$ goes from $w_i$ to $w_j$. 
Also, assume that buds~$\tau_i$ are numbered according to the
ordering of the corresponding base nodes $v_{\tau_i}$ in $W$.
Let these base nodes separate the sequence $w_1, \ldots, w_n$
into parts $Z_1, \ldots, Z_{k+1}$ (some of them may be empty).
In other words, let $Z_i$ be the sequences of nodes from $S$ such that
$w_1, \ldots, w_n = Z_1, v_{\tau_1}, Z_2, \ldots, Z_k, v_{\tau_k}, Z_{k+1}$.
Additionally, put
$A_i := (Z_1 \cup Z_1') \cup V_{\tau_1} \cup \ldots \cup V_{\tau_{i-1}} \cup (Z_i \cup Z_i')$.
Obviously, sets~$A_i$ are self-symmetric, $A_{k+1} = V_G$. The graph~$G[A_1]$ is strongly acyclic
(this readily follows from \refth{strong_acycl_decomp} by putting $Z := Z_1$).
One can show that for each $i \ge 2$ the triple $(A_{i-1}, V_{\tau_{i-1}}, Z_i)$
is a weak separator for $G[A_i]$. Using known decompositions
of $G[V_{\tau_i}]$ these separators can be combined into
a decomposition of $G$. An example is depicted in \reffig{barrier_to_decomp}.

\medskip

Buds that are trimmed by the algorithm are identified in a current
graph but can also be regarded as buds in the original graph~$G$. Namely,
let $H$ be a current graph and $\tau$ be a bud in $H$.
One can see that $(\bar V_\tau, \bar a_\tau)$,
where $\bar a_\tau$ (resp.~$\bar V_\tau$) is the preimage of
$a_\tau$ (resp.~$V_\tau$), is a bud in $G$.
This bud will be denoted by~$\bar \tau$.

Observe that the node sets of preimages $\bar\tau$ of buds~$\tau$ trimmed by
\textsc{Acyclicity-Test} are distinct sets forming a laminar family in $V_G$.
At any moment the current graph~$H$ can be obtained from~$G$ by
trimming the set of inclusion-wise maximal buds (which were discovered up
to that moment). For each such bud~$\bar\tau$ we maintain
an acyclic $v_{\bar\tau}'$-barrier~$\calB^{\bar\tau}$ with the empty $M$-part.

Nodes in $H$ can be of two possible kinds: \emph{simple} and \emph{complex}.
Simple nodes are nodes that were not touched by trimmings, that is, they do not belong
to any of $V_{\bar\tau}$ sets for all trimmed buds $\tau$.
Complex nodes are base and antibase nodes of maximal trimmed buds.

In \refsec{proofs} we prove the following key properties of
\textsc{Acyclicity-Test}:
\begin{enumerate}
	\item[(A)] The (standard directed) subgraph induced by the set of black nodes is acyclic.
	\item[(B)] No arc goes from black node to gray, white or antiblack node.
\end{enumerate}

\textsc{Decompose} consists of two phases: \emph{traversal} and \emph{postprocessing}.
During the first phase we invoke \textsc{Acyclicity-Test} modified as follows.
Suppose the algorithm trims a bud~$\tau$ in~$H$.
First, suppose that the node~$v_\tau$ was simple prior to that trimming.
We construct $\calB^{\bar\tau}$ as follows. Let $B$ be the set of black
simple nodes in $V_\tau$, $\bar\tau_1, \ldots, \bar\tau_k$ be
the preimages (in $G$) of trimmed buds corresponding to base nodes in $V_\tau$.
We argue that putting
$\calB^{\bar\tau} := (B \cup \set{v_\tau'}, \emptyset; \bar\tau_1, \ldots, \bar\tau_k)$
we obtain a required acyclic barrier for $\tau$.
Indeed, property (i) holds by induction; property (ii) follows
from acyclicity of the (standard directed) graph induced by $B$
and the fact that the node~$v_\tau$ (resp. $v_\tau$) cannot
have incoming (resp. outgoing) arcs other than $a_\tau$ (resp. $a_\tau'$).

Situation gets more involved when $v_\tau$ is a complex node
(hence, the algorithm performs several trimmings at this node).
Define $B$ as above. Let $\bar\phi$ be the already trimmed
inclusion-wise maximal bud at $v_\tau$. Consider a
barrier $\calB^{\bar\phi} = (Q, \emptyset; \bar\phi_1, \ldots, \bar\phi_l)$.
We put $\calB^{\bar\tau} := (Q \cup B; \emptyset; \bar\phi_1, \ldots, \bar{\strut\phi_l},
\bar{\strut\tau_1}, \ldots, \bar{\strut\tau_k})$
and argue that a required acyclic barrier for $\bar\tau$ is ready.

First, no arc can go from $Q$ to $B'$. Suppose toward contradiction that
such arc~$uv$ exists. Consider a moment immediately preceding the trimming of~$\phi$.
The tail~$u$ differs from $v_\tau'$ (since $uv \ne a_\tau'$
and $a_\tau'$ is the only arc leaving $v_\tau'$), thus $u$ was black
at that moment. Property~(B) implies that $v$ was either black or antigray.
But this is a contradiction since $v$ is antiblack when $\tau$ is trimmed
and the sets of gray and antigray nodes could not have changed between these
two trimmings.

Second, no arc in $G$ connects $V_{\bar\phi_i}$ and $V_{\bar\tau_j}$.
To see this, suppose toward contraction that $uv$ is such arc, clearly
$uv$ is not the base or antibase arc of $\bar\phi_i$ or $\bar\tau_j$.
Two cases are possible depending on what node, $v_{\bar\phi_i}$~or~$v_{\bar\tau_j}$,
was made black first. Suppose $v_{\bar\phi_i}$ is made black
before $v_{\bar\tau_j}$ (the other case is analogous).
Consider the moment when $v_{\bar\phi_i}$ has just been declared black.
The arc $uv$ in $G$ corresponds to the arc $v_{\bar\phi_i}w$
for some~$w$. The node~$w$ is either black or antigray (due to property~(B)).
In the former case $w = v$; in the latter case $w$ is the antibase node of some already trimmed
bud contained in $\bar\tau_j$ (possibly $\bar\tau_j$ itself). One can see that in
both cases $v_{\bar\tau_j}$ is gray ($v_{\bar\phi_i}$ is an ancestor of $v_{\bar\tau_j}$
in~$F$). Now consider the moment when the algorithm
is about to make $v_{\bar\tau_j}$ black.
The arc~$v'u'$ corresponds to the gray-to-antiblack arc $v_{\bar\tau_j}v_{\bar\phi_i}'$
and hence another trimming at $v_{\bar\tau_j}$ is required. But this is a contradiction as
buds~$\bar\phi_i$ and $\bar\tau_j$ are node-disjoint.

A similar reasoning shows that no arc can go from $Q$ to $V_{\bar\tau_i}$
or from $B$ to $V_{\bar\phi_i}$.

\medskip

When traversal of $G$ is complete the algorithm
builds a final acyclic barrier in $G$. Observe that at that moment
all nodes are black or antiblack. The set of simple black
nodes~$B^*$ in the final graph and the inclusion-wise
maximal trimmed buds $\bar\tau_1^*, \ldots, \bar\tau_k^*$ induce
the acyclic barrier $\calB^* := (B^*, \emptyset; \bar\tau_1^*, \ldots, \bar\tau_k^*)$
in $G$. During the postprocessing phase the algorithm constructs
the desired decomposition of $G$ from acyclic barriers recursively as indicated above.

\begin{remark}
	\textsc{Decompose} procedure relies
	on \textsc{Acyclicity-Test} algorithm and hence the input graph~$G$
	should obey the degree property. One possible workaround
	is to preprocess $G$ as described in \refrem{degree_reduction}.
	What we get after invoking \textsc{Decompose} is a weak
	decomposition~$D$ for such preprocessed graph. A decomposition for $G$ can
	be easily extracted from~$D$. This transformation is straightforward
	and we omit details here.
\end{remark}
	
The presented algorithm also yields a constructive proof to the existence
of a weak separator in a weakly acyclic skew-symmetric graph.
There is, however, a much simpler direct proof based on the same ideas
(it is given in \refsec{proofs}).

\section{Proofs of Separator and Decomposition Theorems}
\label{sec:proofs}

\begin{proofof}{\refth{strong_acycl_decomp}}
	Suppose $G$ has no cycles. Thus $G$ admits topological ordering of nodes (see~\cite{CLR-90}):
	one may assign distinct labels $\pi \colon V_G \to \R$ to the nodes of $G$ so that
	$\pi(u) < \pi(v)$ for every arc $uv \in A_G$. Put $\tilde \pi(v) := \pi(v) - \pi(v')$
	for each~$v \in V_G$. Labeling~$\tilde \pi$ is nowhere zero (since all labels~$\pi$ are distinct)
	and antisymmetric ($\tilde \pi(v) = -\tilde \pi(v')$ for all~$v$). Moreover,
	the skew symmetry of $G$ implies that $\tilde \pi(u) < \tilde \pi(v)$ for all arcs $uv \in A_G$.
	Now consider the set $Z := \setst{v}{\pi(v) > 0}$. Clearly $V_G = Z \sqcup Z'$.
	The induced (standard directed) subgraphs $G[Z]$ and $G[Z']$ are acyclic,
	and no arc leaves~$Z$, as required.

	Conversely, the existence of a partition $Z \sqcup Z'$ of $V_G$ with
	$\deltaout(Z) = \emptyset$ implies that every cycle in $G$
	is contained in $G[Z]$ or $G[Z']$. Since these subgraphs are acyclic,
	theorem follows.
\end{proofof}

We shall need the following result concerning the so-called \emph{regular reachability problem}
(see \cite{GK-96} for more details).
Let $G$ be a skew-symmetric graph with a designated node $s$.
The problem is to find a regular $s$--$s'$ path in $G$ or establish that no such path exists.
One can easily show that if $G$ has an $s$-barrier~$\calB = (S, M; \tau_1, \ldots, \tau_k)$,
then $s'$ is not reachable from $s$ by a regular path.
Indeed, we start from $s$ and need to leave~$S$ in order to reach $s'$.
But after leaving $S$ (via an arc~$a$) we get into one of $V_{\tau_i}$ and can only leave it
by going back to $S$ (because of regularity we are forbidden to use the unique
$V_{\tau_i}$--$S'$ arc~$a'$). Hence, we can never reach $s'$.

Interestingly, the converse statement also holds:

\begin{theorem}[Regular Reachability Criterion, \cite{GK-96}]
\label{th:regular_reach_criterion}
	There exists a regular $s$--$s'$ path in a skew-symmetric graph~$G$
	iff there is no $s$--$s'$ barrier in $G$.
\end{theorem}

We apply this result to prove separator theorems.
Firstly, we need an additional statement:

\begin{lemma}
\label{lm:acycl_reachibility}
	Let $G$ be a weakly acyclic skew-symmetric graph and $s$ be an arbitrary node of $G$.
	Then either $s'$ is not reachable from $s$ by a regular path, or
	$s$ is not reachable from $s'$ by a regular path.
\end{lemma}
\begin{proof}
	For sake of contradiction, suppose that $P$ and $Q$ are regular $s$--$s'$ and
	$s'$--$s$ paths respectively. Then, $P$ cannot be arc-disjoint from both $Q$ and $Q'$
	(since otherwise $P \circ Q$ is a regular cycle).
	Consider the longest prefix of $P$ that is arc-disjoint from $Q, Q'$;
	denote this prefix by $P_0$. Put $P = P_0 \circ P_1$ and let $a$ be the first arc of $P_1$.
	This arc is either contained in $Q$ or $Q'$ but not in both (as $Q$ is regular).
	Assume, without loss of generality, that $a$ belongs to $Q$ and let
	$Q_1$ be the suffix of $Q$ starting with~$a$ (that is, $a$ is the first arc of~$Q_1$). 
	Combine $P_0$ and $Q_1$ into the cycle $C := P_0 \circ Q_1$. This cycle is regular
	since its initial part $P_0$ is regular and arc-disjoint from both $Q_0$ and $Q_0'$.
	The regularity of $C$ contradicts the weak acyclicity of $G$, and the claim follows.
\end{proof}

\begin{proofof}{\refth{weak_separator}}
	The proof is by induction on $\abs{V_G}$.
	
	Applying \reflm{acycl_reachibility} we get a node $s$ such that $s'$ is not reachable
	from $s$ by a regular path. Let $\calB = (S, M; \tau_1, \ldots, \tau_k)$
	be an $s$--$s'$ barrier (which exists due to \refth{regular_reach_criterion}).
	Two cases are possible:
	
	\emph{Case 1: $k = 0$.}
	Then we apply the induction hypothesis to $G[M]$.
	If $G[M]$ is strongly acyclic, then so is $G$. (Moreover, if $(Z, Z')$
	is partition of $G[M]$ as in \refth{strong_acycl_decomp} then
	$(Z \cup S, Z' \cup S')$ is a similar partition for $G$.)
	Otherwise $G[M]$ is weakly acyclic; consider a weak separator $(A, B, Z)$
	for $G[M]$. Now $(A, B, Z \cup S)$ is a weak separator for $G$
	and the induction follows.
	
	\emph{Case 2: $k > 0$.}
	Let $\tilde G$ denote the graph obtained from $G$ by trimming all buds $\tau_i$
	(it is clear that the order of trimmings is unimportant).
	Let $W$ be the set of nodes formed by adding nodes~$v_{\tau_i}$ to~$S$.
	We argue that the (standard directed) subgraph $\tilde G[W]$ is acyclic.
	Indeed, suppose $C$ is a circuit in $\tilde G[W]$. Since $W \cap W' = \emptyset$
	this circuit is regular (in $\tilde G$) and applying the restoration procedure 
	from~\refsec{algorithms} one can transform $C$ into a regular circuit in $G$~--- a contradiction.
	
	Since $\tilde G[W]$ is acyclic we may consider its topological ordering
	$w_1, \ldots, w_n$ (no arc goes from $w_i$ to $w_j$ for $i > j$).
	Let $w_j$ be the node with the largest index that is the base node of some bud~$\tau$ of $\calB$.
	Let $D$ be the set of nodes in $S$ preceding $w_j$ (with respect to the topological order),
	let $Z$ be the set of nodes in $W$ following $w_j$. Put $B := V_\tau$, $A := V \setminus (B \cup Z \cup Z')$.
	
	We claim that $(A, B, Z)$ is a weak separator for $G$.
	
	To see this, we check the properties of weak separator one by one.
	Clearly $A = A'$, $B = B'$, and sets $A, B, Z, Z'$ form a partition of $V_G$;
	$w_j$ is the base node of $\tau$ and hence the arc $a_\tau$ goes from
	$A$ to $B$ (the symmetric arc $a_\tau'$ goes from $B$ to $A$).
	Thus, sets $A$, $B$ are connected by a pair of arcs.
	No other pair of arcs can go between $A$ and $B$: if $a \ne a_\tau$ is an
	arc from $A$ to $B$ then it either connects $V_{\tau_j}$ with some $V_{\tau_i}$
	or goes from $S$ to $V_\tau$ or goes from $D'$ to $V_\tau$.
	The first two cases are forbidden by the definition of barrier.
	Consider the latter case. In graph~$G$ the symmetric arc $a'$ goes
	from $V_\tau$ to $D$. Hence, the image of $a'$ in $\tilde G$
	goes from $w_j$ to $D$ violating the topological ordering of $W$.
	
	Now we prove that no arc can leave~$Z$. Suppose, for the sake of contradiction, that $a$
	is such arc. It is clear that $Z \subseteq S$.
	The arc~$a$ cannot go to $M \cup S'$ or one of $V_{\tau_i}$ (by definition of barrier), thus
	it should go to $S \setminus Z$. The image of $a$ in $\tilde G$, hence, should go to $D$.
	This again is a contradiction.
\end{proofof}

\begin{proofof}{\refth{strong_separator}}
	Sufficiency being obvious, we show necessity.
	Let $G$ be a strongly connected weakly acyclic skew-symmetric graph.
	Consider an arbitrary weak separator $(A, B, Z)$. It follows that $Z = \emptyset$
	(since a strongly connected graph cannot have a nontrivial directed cut).
	Let $a$ (resp. $b$) be the head of the (unique) arc going from $B$ to $A$
	(resp. from $A$ to $B$). Consider arbitrary nodes $u \in A$, $v \in B$.
	Since $G$ is strongly connected, one should have a cycle containing both $u$ and $v$.
	Restricting this cycle on $G[A]$ and $G[B]$ we get $a$--$a'$ and $b$--$b'$ paths
	containing nodes $u$ and $v$ respectively. Hence, $G[A]$ is $a$-connected
	and $G[B]$ is $b$-connected.
\end{proofof}

\begin{proofof}{\refth{acycl_decomp}}
	Consider a partition of $V_G$ into strongly connected components $Q_1, \ldots, Q_k$ and the
	(standard directed) component graph $G_C$ formed by removing arcs inside
	components and contracting components~$Q_i$
	into composite nodes.  Each component $Q_i$ is either self-symmetric ($Q_i' = Q_i$)
	or \emph{regular} ($Q_i \cap Q_i = \emptyset$). In the latter case there exists~$j$
	such that $Q_i' = Q_j$. Hence, the collection of regular components can be partitioned into
	pairs of symmetric ones. 

	Since $G$ is weakly symmetric, each regular component~$Q_i$ is trivial (consists of
	a single node of $G$). Let $W$ denote the set of nodes in regular components
	and $\set{B_1, \ldots, B_n}$ be the collection of all self-symmetric components.
	As $G_C$ is acyclic one may construct a topological labeling of components, that is,
	assign distinct labels $\pi$ to the components satisfying $\pi(K) < \pi(L)$
	for each arc going from $K$ to $L$.
	Note, that for distinct $i, j$ component~$B_i$ is not reachable from $B_j$ in $G_C$
	(if $P$ is an $B_i$--$B_j$ path then $P'$ is an $B_j$--$B_i$ path, contradicting the
	acyclicity of $G_C$). Therefore, we may assign labels $\pi$ to self-symmetric components
	in an arbitrary way. We set all these labels to zero (while the labels of regular components are
	assumed to be distinct).

	Now we transform $\pi$ so as to make it antisymmetric. To this aim
	we define $\tilde \pi(C) := \pi(C) - \pi(C')$ for each component~$C$ (both regular and self-symmetric).
	Clearly, $\tilde \pi$ is antisymmetric and the new labels of symmetric components are still zero
	(while other labels are nonzero). 
	For every arc going from component~$K$ to component~$L$ there is a symmetric arc from $L'$ to $K'$
	and one has $\pi(K) < \pi(L)$, $\pi(L') < \pi(K')$ thus proving $\tilde \pi(K) < \tilde \pi(L)$.
	Hence, $\tilde \pi$ is again a topological labeling on $G_C$. 
	Define $Z$ to be the set of nodes in regular components with positive label~$\tilde \pi$.
	Theorem now follows from the properties of $\tilde \pi$.
\end{proofof}

\section{Correctness of Acyclicity Test}
\label{sec:correctness}

We simultaneously prove (using induction on the number of steps
performed by the algorithm) the following properties:
\begin{enumerate}
	\item[(A)] The (standard directed) subgraph induced by the set of black nodes is acyclic.
	\item[(B)] No arc goes from black node to gray, white or antiblack node.
	\item[(C)] Each time an arc from gray node~$u$ to antiblack node $v$ is discovered $v'$ is an ancestor of $u$.
	\item[(D)] Bud trimming preserves ancestors in $F$ (that is, if $u$ is an ancestor of $v$
	in $F$ before trimming and trimming does not remove neither $u$ nor $v$, $u$ will be an ancestor of $v$
	after trimming).
\end{enumerate}

To prove property~(C) consider a current graph~$H_1$
and suppose toward contradiction that $u$ is gray, $v$ is antiblack, there exists
an arc~$uv$ in $H_1$ but $v'$ is not an ancestor of $u$.
Hence $v'$ was marked as black before $u$ has been discovered. Consider the moment when
$v'$ was made black and the corresponding current graph~$H_0$. The arc~$v'u'$ exists in~$H_0$
(no trimming applied by the algorithm while going from $H_0$ to $H_1$ could
affect $u$ or $v$). But $v'u'$ is a black-to-white arc at the considered moment contradicting to property~(B).
Hence, all trimmings performed by the algorithm operate with well-defined buds.

Property (D) follows trivially.

Consider property~(B). Clearly it is maintained while algorithm changes the colors of nodes without
performing trimmings. Consider a trimming performed at a base node~$u$. This trimming
does not add new black nodes. If $a$ is an arc with black tail it is either left unchanged by the trimming
or redirected toward~$u'$. In the latter case it becomes black-to-antigray and property~(B) holds.

Finally we prove property~(A). Let $B$ be the set of black nodes. Let each node~$v \in B$ be
assigned a moment of time~$f(v)$ when it has become black. (As in standard DFS,
these moments are just arbitrary increasing integers.) We claim that
numbers $f(v)$ give the topological ordering of $B$, that is, $f(u) > f(v)$ for each arc~$uv$.
Observe that when a new node~$v$ becomes black it gets a label that is larger than
all other existing labels. 
Thus, the property is satisfied for all arcs leaving~$v$. Consider an arbitrary incoming arc~$uv$.
The node~$u$ cannot be black: otherwise prior to making~$v$ black the algorithm had
a black-to-gray arc~$uv$, which is forbidden by property~(B).
property~(A) remains valid when the algorithm trims an arbitrary bud~$\tau$
since no new black-to-black arcs are created.

We are now ready to prove the correctness \textsc{Acyclicity-Test}:

\begin{lemma}
\label{lm:cycles_preserved}
	If a current graph of the algorithm has a regular circuit before
	trimming it still has one after trimming.
\end{lemma}
\begin{proof}
	Suppose the algorithm trims a bud~$\tau$ in a current graph~$H$.
	Consider a moment immediately preceding such trimming.
	Let $C$ be a regular circuit in $H$. In case $C$ does not
	intersect $V_\tau$ it is obviously preserved under trimming.
	
	Otherwise, without loss of generality one may assume that $C$
	contains $v_\tau'$ or a black node from~$V_\tau$.
	We argue that $C$ contains an antigray node.
	If $C$ contains $v_\tau'$, then we are done.
	Otherwise, let us go along~$C$ and examine the colors of nodes.
	All nodes of $C$ cannot be black (due to property~(A)).
	Hence, from property~(B) it follows that at some point we have
	a black-to-antigray transition, as required. 
	
	Let $x$ be the antigray node of $C$ such that the unique $v_\tau'$--$x$ path~$R$
	in $F$ is as short as possible. If $C$ does not intersect $R'$, then we are done.
	Otherwise, replace $C$ by $C'$ and try again (the length of $R$ decreases
	since $C$ is regular and hence, due to degree property, cannot contain a pair
	of symmetric nodes $x, x'$). Finally we get a regular circuit~$C$, an
	antigray node~$x$, and a $v_\tau'$--$x$ path~$R$ such that $C$
	contains $x$ and $C$ does not intersect $R'$.

	We go from $x$ along $C$ until reaching $V_\tau$ at $y$.
	Let $Q$ be a regular $y$--$v_\tau'$ path in $G[V_\tau]$;
	let $P$ be the $x$--$y$ segment of $C$. Combine these three paths
	together by putting $K := P \circ Q \circ R$.
	$K$ is a regular cycle in $H$ that is preserved under trimming. Hence,
	the trimmed graph contains a regular circuit, as required.
\end{proof}	

\begin{theorem}
	A skew symmetric graph~$H$ is weakly acyclic iff \textsc{Acyclicity-Test}
	reports no regular cycle in it.
\end{theorem}
\begin{proof}
	The necessity is straightforward and
	has been already established in \refsec{algorithms}.
	To prove sufficiency assume toward contradiction that $G$ has a regular circuit
	but the algorithm did not discover one. By \reflm{cycles_preserved}
	the presence of regular circuits in preserved during the course of the algorithm.
	In the resulting graph all nodes are either black or antiblack.
	Properties~A and B imply that this graph is strongly acyclic (see~\refth{strong_acycl_decomp}).
	Hence, one gets a contradiction.
\end{proof}

\section{Efficient Implementations}
\label{sec:impl}

To implement \textsc{Acyclicity-Test} efficiently we
borrow some standard techniques from \cite{GKT-99,GK-96}. We assume that
it is possible to obtain a mate for any given node or arc in $O(1)$ time.
All graphs are assumed to be represented by adjacency lists.
That is, for any node~$v$ all arcs leaving $v$
are organized in a double-linked list attached to $v$. These lists allow to enumerate
all arcs leaving a given node in time proportional to its out-degree.
Note that we do not maintain lists of incoming arcs explicitly. Instead,
to enumerate the incoming arcs of $v$ we enumerate arcs leaving $v'$ and
apply symmetry.

Let $H$ be a current graph at some point of execution of \textsc{Acyclicity-Test}.
For a node $x$ in $G$ let $\hat x$ be the node in $H$ defined as follows:
if $x$ is a simple node in $H$, then $\hat x := x$; otherwise $\hat x := v_{\bar\tau}$
where $\bar\tau$ is the maximal trimmed bud in $G$ containing $x$.
The arcs of $H$ are represented by their preimages in~$G$.
More formally, consider an arbitrary arc~$a$ in $H$ and let $uv$ be the
corresponding arc (preimage) in $G$. Trimmings could have changed the head and the
tail of $uv$. One can easily check that the arc~$a$ in $H$ goes from $\hat u$ to $(\hat v)'$.
To compute $\hat x$ by $x$ efficiently we use an instance of \emph{disjoint set union}
data structure (see \cite{CLR-90}) and denote it by $\calF$. These disjoint sets
are the node sets of maximal trimmed buds in $G$. An argument
as in \cite{GK-96, GKT-99} shows that operations performed by our algorithm on $\calF$
fall into a special case admitting $O(1)$ cost for \emph{unite} and
\emph{find} calls (this implementation is given in~\cite{GT-85}).

Next we consider bud trimming operation and discuss its implementation.
Consider a trimming of a bud~$\tau$ in a current graph~$H$.
Firstly, we update $\calF$ by performing \emph{unite} on it to
reflect the changes in the structure of maximal buds.
Secondly, we need to update the graph adjacency lists.
The naive approach would be as follows.
Enumerate all arcs leaving $v_\tau$. Construct a new list of outgoing arcs (skipping
arcs in $\gamma(V_\tau)$) and attach it to $v_\tau$ (replacing the old list).
This approach, however, is inefficient since an arc can be scanned many times
during the execution of the algorithm. To do better, we concatenate the lists of
outgoing arcs for the nodes in $V_\tau$ and attach the resulting list to $v_\tau$.
This takes $O(\abs{V_\tau})$ time. Unfortunately, this also yields an additional issue:
the arcs in $\gamma(V_\tau)$ (which are normally removed by trimming) remain
in the current graph. We cannot identify such arcs during trimming procedure since that
would require to scan all arcs leaving $v_\tau$ and take too much time.
Instead, we use \emph{lazy deletion} strategy as in \cite{GK-96}: call an arc of the
initial graph $G$ \emph{dead} if it is contained in $\gamma(V_{\bar\tau})$ for some
maximum trimmed bud~$\bar\tau$. We admit the presence of dead arcs in our lists
but remove such arcs as soon as we discover them. Dead arcs~$uv$ can easily be detected
by checking if $\hat u = \hat v$. Since an arc can be removed at most
once, the running time is not affected.

To maintain the forest~$F$ we store for each non-root node~$v$ of $F$ the unique forest
arc~$q(v)$ entering~$v$. This information about $F$ allows
to construct the node set~$V_\tau$ of a discovered bud~$\tau$ 
in $O(\abs{V_\tau})$ time. Suppose that on the current step
the algorithm examines an arc~$uv$ in a current graph~$H$ and finds out that $v$ is antiblack.
Then, $v'$ is an ancestor of $u$ (as shown in \refsec{proofs}). We trace the
corresponding $u$--$v'$ path in $F$ in backward direction hence obtaining $V_\tau$.
No additional processing is required to update $F$ when a bud~$\tau$ is trimmed:
the nodes in $V_\tau \setminus \set{v_\tau, v_\tau'}$ vanish from the current graph
and these values~$q(v)$ are no longer used.

Let $n$ (resp. $m$) denote the number of nodes (resp. arcs) in the original graph.
The described implementation scans each arc at most once and takes $O(1)$ time
for each examination. Also, $O(n)$ additional time units are required for initialization
and other auxiliary actions. The time consumed by trimmings is proportional to
$O(\sum_i \abs{V_{\tau_i}})$ where $\tau_i$ are all trimmed buds. This sums telescopes to $O(n)$.
The restoration procedure performed on a regular cycle by our algorithm has
the running time $O(n)$ and is essentially the same as the corresponding
one in \cite{GK-96}, so we omit details here. Finally, we conclude that \textsc{Acyclicity-Test}
can be implemented to run in $O(m + n)$ time.

\medskip

Now let us estimate the complexity of \textsc{Decompose}. In order to be efficient
we need a compact way of storing decompositions. (For example, listing all the corresponding
sets $A$, $B$, $Z$ explicitly may require $\Theta(n^2)$ space.)
To do this, we only store $Z$-part for each node of decomposition tree. Since these
sets are disjoint, the linear bound on the size of decomposition follows. Obviously,
we may still obtain $A$- and $B$-parts (if required) of any node~$x$ in decomposition tree by traversing the
two subtrees rooted at children of $x$ and uniting the corresponding $Z$- and $Z'$-parts.
Each $Z$-part stored by the algorithm is organized as a double-linked list.
During the postprocessing phase the algorithm converts acyclic barriers into weak decomposition of $G$.
To avoid invoking topological sort of each bud we collect values~$f(v)$
for all nodes~$v$ that become black during traversal phase (as described in~\refsec{proofs}).
These time labels induce topological order (as in standard DFS algorithm).
Hence, to postprocess a bud $\tau$ the algorithm requires $O(\abs{V_\tau})$ time units.
Therefore, the postprocessing phase runs in $O(n)$ and the total complexity
of \textsc{Decompose} is linear.

\section{Concluding Remarks}

We have studied the structure of weakly acyclic bidirected and skew-symmetric graphs.
The obtained decomposition theorems combine the notions of topological
ordering (as in case of standard directed graphs)
and barrier (which is a standard tool for working with regular reachability
problems in skew-symmetric graphs). We have adopted the algorithm of Gabow,
Tarjan, and Kaplan to test weak acyclicity in linear time. Moreover,
we have proposed a variation of this method to build weak acyclic
decomposition in linear time. The problem of extending such algorithm
to the case of strong decomposition remains open.

\section{Acknowledgments}
The author is thankful to Alexander Karzanov for constant attention,
collaboration, and many insightful discussions.

\nocite{*}
\bibliographystyle{plain}
\bibliography{main}

\end{document}